\documentclass[12pt]{article}
\usepackage{amssymb}
\usepackage{amsmath}
\usepackage{amsbsy}
\usepackage{amsthm}
\usepackage{rotating}

\newcommand{\be}{\begin{equation}}
\newcommand{\ee}{\end{equation}}
\newcommand{\ba}{\begin{eqnarray}}
\newcommand{\ea}{\end{eqnarray}}
\newcommand{\ban}{\begin{eqnarray*}}
\newcommand{\ean}{\end{eqnarray*}}

\newtheorem{theo}{Theorem}[section]
\newtheorem{prop}[theo]{Proposition}

\newtheorem{cor}[theo]{Corollary}

\begin{document}
\title{Hyperbolic groups admit proper affine isometric actions on $l^p$-spaces}
\author{Guoliang Yu\footnote{Partially supported by NSF and NSFC.}}
\date{ }
\maketitle

\section{Introduction}

Let $X$ be a Banach space and $\Gamma$ be a countable discrete
group. An affine and isometric action $\alpha$ of $\Gamma$ on $X$
is said to be proper if $lim_{g\rightarrow \infty} \| \alpha(g)
\xi \| =\infty$ for every $\xi\in X$. If $\Gamma$ admits a proper
isometric affine action on Hilbert space, then $\Gamma$ is said to
be of Haagerup property [9] or a-T-menable [12].

Bekka,  Cherix and  Valette proved that an amenable group admits a
proper affine isometric action on Hilbert space [3]. This result
has important applications to K-theory of group $C^*$-algebras
[13] [14].

It is well known that an infinite Property (T) group doesn't admit
a proper affine isometric action on Hilbert space. The purpose of
this paper is to prove the following result.

\begin{theo}
If $\Gamma$ is  a hyperbolic group, then there exists $2\leq
p<\infty$ such that $\Gamma$ admits a proper affine isometric
action on an $l^p$-space.
\end{theo}

We remark that the constant $p$ depends on  the hyperbolic group
$\Gamma$ (in the special case that $\Gamma$ is the fundamental
group of a negatively curved compact manifold, $p$ depends on the
dimension of the manifold), and $p$ is strictly greater than $2$
if the hyperbolic $\Gamma$ is infinite and has Property (T).
Recall that a theorem of A. Zuk states that hyperbolic groups are
generically of Property (T) [22].

In [1], Bader and Gelander studied Property (T) for $L^p$-spaces.
 Their work has extremely interesting applications in Fisher and Margulis' theory
 of
local rigidity [6]. Bader and Gelander raised the question if any
affine isometric action of a Property (T) group on an $L^p$-space
has a fixed point (Question 12 in [1]). Theorem 1.1 implies that
the answer to this question is negative for infinite hyperbolic
groups with Property (T).

The proof of Theorem 1.1 is based on a construction of Igor
Mineyev [18] and is reminiscent of Alain Connes' construction of
Chern character of finitely summable Fredholm modules for rank one
groups [5].

The author wishes to  thank Igor Mineyev for very helpful comments
on the exposition of this note,  Erik Guentner for bringing [1] to
my attention, and Nigel Higson for pointing out that an
unpublished result Y. Shalom implies that $Sp(n,1)$ admits a
proper affine isometric action on some uniformly convex Banach
space.

\section{Hyperbolic groups and bicombings.}

In this section, we recall the concepts of hyperbolic groups and
bicombings.

\subsection{Hyperbolic groups.}

Let $\Gamma$ be a finitely generated group.
 Let $S$ be a  finite generating
set for $\Gamma$. Recall that the Cayley graph  of $\Gamma$ with
respect to $S$ is the graph $G$  satisfying the following
conditions:
\begin{itemize}
\item [(1)] the set of vertices in $G$, denoted by $G^{(0)}$, is $\Gamma$;
\item [(2)] the set of edges is $\Gamma \times S$, where each edge
  $(g,s)\in \Gamma\times S$ spans the vertices $g$ and $g s.$
\end{itemize}

We endow  $G$  with the path metric $d$ induced by assigning
length 1 to each edge. Notice that $\Gamma$ acts freely,
isometrically and cocompactly on $G$.
 A geodesic path in $G$ is a shortest edge path.
The restriction of the path metric $d$ to $\Gamma$ is called the
word metric.

A finitely generated group $\Gamma$ is called  hyperbolic, if
there exists a constant $\delta\geq 0$ such that all the geodesic
triangles in $G$ are $\delta$-fine in the following sense: if $a$,
$b$, and $c$ are vertices in $G$, $[a,b]$, $[b,c]$, and $[c,a]$
are geodesics from $a$ to $b$, from $b$ to $c$, and from $c$ to
$a$, respectively, and points $\bar{a}\in [b,c]$, $v,\bar{c}\in
[a,b]$, $w,\bar{b}\in [a,c]$ satisfy
$$d(b,\bar{c})=d(b,\bar{a}),\quad d(c,\bar{a})=d(c,\bar{b}),\quad
d(a,v)=d(a,w)\leq d(a,\bar{c})=d(a,\bar{b}),$$ then $d(v,w)\leq
\delta$.

The above definition of hyperbolicity  does not depend on the
choice of the  finite generating set $S$.
 See [8]
 for other equivalent definitions.

For vertices $a$, $b$, and $c$ in $G$, the Gromov product is
defined by
$$(b|c)_a := d(a,\bar{b})=  d(a,\bar{c})=
  \frac{1}{2}\Big[d(a,b)+d(a,c)-d(b,c)\Big].$$
The Gromov product can be used to
 measure the degree of cancellation in the multiplication
of group elements in $G$.

\subsection{Bicombings.}
Let $\Gamma$ be a finitely generated group. Let $G$ be its Cayley
graph with respect to a finite generating set. A  bicombing $q$ in
$G$ is a function assigning to each ordered pair $(a,b)$ of
vertices in $G$ an oriented edge-path $q[a,b]$ from $a$ to $b$. A
bicombing $q$ is called  geodesic, if each path $q[a,b]$ is
geodesic, i.e. a shortest edge path. A bicombing $q$ is
$\Gamma$-equivariant if $q[g\cdot a, g\cdot b]= g\cdot q[a,b]$ for
each $a,b\in G^{(0)}$ and each $g\in \Gamma$.

\section{A construction of Mineyev.}

The purpose of this section is to recall Mineyev's contruction for
hyperbolic groups  and its properties [18].

Let $\Gamma$ be a hyperbolic group and $G$ be a Cayley graph of
$\Gamma$ with respect to a finite  generating set.
 We endow $G$ with  the  path metric $d$,
 and identify $\Gamma$ with $G^{(0)}$,  the set of vertices of $\Gamma$.
Let $\delta\ge 1$ be a positive integer such that all the geodesic
triangles in $G$ are $\delta$-fine.

The  ball $B(x,R)$ is the set of all vertices at distance at most
$R$ from the vertex $x$. The sphere $S(x,R)$ is the set of all
vertices at distance $R$ from the vertex $x$. Pick an equivariant
geodesic bicombing~$q$ in~$G$.  By $q[a,b](t)$ we denote the point
on the geodesic path $q[a,b]$ at distance $t$ from~$a$. Recall
that $C_0 (\Gamma,\mathbb{Q})$ is the space of all finitely
supported
 $0$-chains (in $\Gamma=G^{(0)}$)  with coefficients in
$\mathbb{Q}$, i.e. $C_0 (\Gamma,\mathbb{Q})= \{\sum_{\gamma\in
\Gamma} c_{\gamma}\gamma: c_{\gamma}\in \mathbb{Q}\}$, where
$\sum_{\gamma\in \Gamma} c_{\gamma}\gamma $ is finitely supported.

For each $p\geq 1$, endow $C_0(\Gamma,\mathbb{Q})$ with the
$l^p$-norm $\|\cdot\|_p$. We identify $\Gamma$ with the standard
basis of $C_0(\Gamma,\mathbb{Q})$. Therefore the left action of
$\Gamma$ on itself induces a left action on~$C_0(G,\mathbb{Q})$.

For $v, w\in \Gamma$, the  flower at $w$ with respect to $v$ is
defined to be
$$ Fl(v, w) := S(v, d(v,w))\cap B(w,\delta)\subseteq \Gamma.$$

For each $a\in \Gamma$, we define $pr_a: \Gamma\rightarrow \Gamma$
by:
\begin{itemize}
\item [(1)] $pr_a (a) := a;$
\item [(2)] if $b\neq a$, $pr_a(b) := q[a,b](t)$, where $t$ is the largest
integral multiple of $10\delta$ which  is strictly less than
$d(a,b)$.
\end{itemize}

Now for each pair $a,b\in \Gamma$, we define a $0$-chain $f(a,b)$
in $\Gamma$ inductively on the distance $d(a,b)$ as follows:
\begin{itemize}
\item [(1)] if $d(a,b)\leq 10\delta$, $f(a,b) :=b$;
\item [(2)] if $d(a,b) >10\delta$ and $d(a, b)$ is not an integral
multiple of $10\delta$, let $f(a,b) := f(a, pr_a (b));$
\item [(3)] if $d(a,b)> 10\delta$ and $d(a,b)$ is an integral
 multiple of $10\delta$,
 let
$$ f(a,b) := \frac{1}{\# Fl(a,b)}\sum_{x\in Fl(a,b)} f(a, pr_a (x)).$$
\end{itemize}

The following result is due to Mineyev [18].

\begin{prop}
The function $f:\Gamma \times \Gamma \to C_0(\Gamma,\mathbb{Q})$
defined above satisfies the following conditions.
\begin{itemize}
\item [(1)] For each $a,b\in \Gamma$, $f(b,a)$ is a  convex combination,
i.e. its coefficients are non-negative and sum up to~1.
\item [(2)] If $d(a,b)\geq 10\delta$, then
$supp\,f(b,a)\subseteq B(q[b,a](10\delta),\delta)\cap
S(b,10\delta)$.
\item [(3)] If $d(a,b)\leq 10\delta$, then $f(b,a)=a$.
\item [(4)] $f$ is $\Gamma$-equivariant, i.e.
$f(g\cdot b, g\cdot a)= g\cdot f(b,a)$ for any $g,a,b\in \Gamma$.
\item [(5)] There exist constants $L\geq 0$ and $0\leq \lambda<1$
such that, for all $a,a',b\in \Gamma$,
$$\|f(b,a)-f(b,a')\|_1\leq L\,\lambda^{(a|a')_b}.$$
\end{itemize}
\end{prop}

Let $p\geq 2$. For each pair $b,a \in \Gamma$, define
$$ h(b,a) =\frac{1}{ \| f(b,a)\|_{p}} f(b,a),$$
where $f$ is as in Proposition 3.1.

\begin{cor}
The function $h:\Gamma \times \Gamma \to C_0(\Gamma,\mathbb{Q})$
defined above satisfies the following conditions.
\begin{itemize}
\item [(1)] For each $a,b\in \Gamma$, $\|h(b,a)\|_{p}=1$.
\item [(2)] If $d(a,b)\geq 10\delta$, then
$supp\,h(b,a)\subseteq B(q[b,a](10\delta),\delta)\cap
S(b,10\delta)$.
\item [(3)] If $d(a,b)\leq 10\delta$, then $h(b,a)=a$.
\item [(4)] $h$ is $\Gamma$-equivariant, i.e.
$h(g\cdot b, g\cdot a)= g\cdot h(b,a)$ for any $g,a,b\in \Gamma$.
\item [(5)] There exist constants $C\geq 0$ and $0\leq \rho<1$
such that, for all $a,a',b\in \Gamma$,
$$\|h(b,a)-h(b,a')\|_{p}\leq C\, \rho^{(a|a')_b}.$$
\end{itemize}
\end{cor}

{\noindent {\it Proof:}}
 (1), (2), (3) and (4) of Corollary 3.2 follow  from
Proposition 3.1.

By (2) of Proposition 3.1, we have
$$\# supp \,h(b,a) \leq \# S(b,10 \delta), \,\, \# supp \,h(b,a') \leq \# S(b,10
\delta).$$ It follows that
$$\|h(b,a)-h(b,a')\|_{p}\leq  \,\,2(\#
S(b,10\delta))^{\frac{1}{p}}\,\,\, \|h(b,a)-h(b,a')\|_{1}.$$ Now
(5) of  Corollary 3.2 follows from (5) of Proposition 3.1.

\qed

\section{Proof of the main result.}

In this section, we prove Theorem 1.1.

{\noindent {\it Proof of Theorem 1.1:}}

Let $\upsilon>0$ such that $\# B(x, r)\leq \upsilon^r$ for all
$x\in \Gamma$ and $r>0$. Let $\rho$ be as in Corollary 3.2. Choose
$p\geq 2$ such that $\rho^p \upsilon < \frac{1}{2}.$

Let $l^p(\Gamma)$ be the completion of $C_0(\Gamma, \mathbb{Q})$
with respect to the norm $\|\cdot \|_p$. Notice that the $\Gamma$
action on  $C_0(\Gamma, \mathbb{Q})$ can be extended to an
isometric action  on $l^p(\Gamma)$.

 Let
$$X=\{\xi:\Gamma\rightarrow l^p(\Gamma): \|\xi\|_p=(\sum_{\gamma \in \Gamma}
\| \xi (\gamma)\|^p )^{\frac{1}{p}}<\infty\}.$$ Observe that $X$
is isometric to $l^p(\Gamma \times \Gamma).$

 Let $\pi$ be the isometric action of
$\Gamma$ on $X$ defined by:

$$(\pi(g)\xi)(\gamma)=g(\xi (g^{-1}\gamma))$$
for all $\xi\in X$ and  $g, \gamma\in \Gamma.$

Define $\eta\in X$ by:
$$\eta(\gamma)= h(\gamma,e)$$
for all $\gamma\in \Gamma$, where $e$ is the identity element in
$\Gamma$.

For each $g\in \Gamma$, by  Corollary 3.2 and the choice of $p$,
we have:
$$\|\pi(g)\eta-\eta\|_{p}^p= \sum_{\gamma\in\Gamma} \| g
(h(g^{-1}\gamma,e))- h(\gamma,e)\|_{p}^p~~~~~~~~~~~~$$
$$\leq \sum_{\gamma\in\Gamma} \|
h(\gamma,g)- h(\gamma, e)\|_{p}^p$$
$$\leq \sum_{\gamma\in\Gamma} C^p\rho^{p (d(\gamma,e)-d(g,e))}~~~~$$
$$\leq \sum_{n=0}^{\infty} C^p \rho^{p (n-d(g,e))}
\upsilon^n~~~~$$ $$\leq \,\,\,2C^p\,\,\, \rho^{-p
d(g,e)}\,\,\,\,\,\,\,\,\,.\,\,\,\,\,\,\,\,\,\,\,\,$$ It follows
that $\pi(g)\eta-\eta$ is an element in $X$ for each $g\in
\Gamma$.

We now define an affine isometric action $\alpha$ on $X$ by
$\Gamma$ by:

$$\alpha(g) \xi=\pi (g) \xi + \pi(g)\eta -\eta$$
for all $\xi \in X$ and $g\in\Gamma$.

If $\gamma$ is a vertex on the oriented geodesic $q[g,e]$
satisfying $d(\gamma,e)\geq 10\delta$ and $d(\gamma, g)\geq
10\delta$, we have
$$B(q[\gamma,e](10\delta), \delta) \cap B(q[\gamma,g](10\delta),
\delta)=\varnothing.$$

 Otherwise if there exists $z \in B(q[\gamma,e](10\delta),
\delta) \cap B(q[\gamma,g](10\delta), \delta),$ then
\begin{eqnarray*}
d(g, e)&\leq &d(g, z)+ d(z,e)\\
& \leq &(d(g, q[\gamma,g](10\delta))+\delta))+
(\delta+d(q[\gamma,e](10\delta),e)\\
&=&((d(g,\gamma)-10\delta)+\delta)+(\delta + (
d(\gamma,e)-10\delta))\\
& =& d(g,e)-18\delta.
\end{eqnarray*}

 This is a contradiction since $\delta>0$.

 By (2)  of Corollary 3.2, we have
$$supp \,\,h(\gamma, g) \cap supp \,\,h(\gamma,e)=\varnothing$$
if $\gamma$ is a vertex on the oriented geodesic $q[g,e]$
satisfying $d(\gamma,e)\geq 10\delta$ and $d(\gamma, g)\geq
10\delta$.

It follows  that there exist at least $d(g,e)-100\delta$ number of
vertices $\gamma$ on the oriented path $q[g, e]$ such that
$$\| g (h(g^{-1}\gamma,e))-
h(\gamma,e)\|_{p}=\|h(\gamma,g)-h(\gamma,e)\|_p\geq 1.$$

Hence
$$\|\pi(g)\eta-\eta\|_p^p\geq d(g,e)-100\delta$$
for all $g\in \Gamma.$

As a consequence,  for every $\xi \in X$, we have
$$\|\alpha (g) \xi-\pi (g) \xi \|_p \rightarrow \infty$$
as $g\rightarrow \infty$.

This, together with the fact that $\pi (g)$ is an isometry,
implies that $\alpha$ is proper.

\qed

We should mention that it remains an open question if
$SL(n,\mathbb{Z})$ admits a proper affine isometric action  on
some uniformly convex Banach space for $n\geq 3$. A positive
answer to this question would have interesting applications to
K-theory of group $C^*$-algebras [16].

\noindent Department of Mathematics

\noindent 1326 Stevenson Center

\noindent Vanderbilt University

\noindent Nashville, TN 37240, USA

\noindent e-mail: gyu@math.vanderbilt.edu

\end{document}